\input amstex
\input amsppt.sty
\documentstyle{amsppt}
\magnification\magstep1

\def\Z{\text{\bf Z}}

\NoBlackBoxes
\hoffset= 0.0in
\voffset= 0.0in
\hsize=32pc
\vsize=43pc

\topmatter 	\title Compact group actions 
that raise dimension to infinity 	
\endtitle
\footnote{Research partially supported by the National Science 
Foundation}
	\author A. N. Dranishnikov and J. E. West 	\endauthor

	\affil 
University of Florida and Cornell University 	\endaffil

        \address Department of Mathematics, University of Florida
 358 Little Hall, PO Box 118105, Gainesville, FL 32611\endaddress        

	\address Department of Mathematics, White Hall, Cornell 
University, 		Ithaca, NY 14853 	\endaddress

	\email dranish\@math.ufl.edu west\@math.cornell.edu 
	\endemail



	\abstract\nofrills THEOREM. {\it For every prime $p$ and each $n=2, 3, ... 
\infty$, 		there is an action of $G=\prod 
_{i=1}^{\infty}(Z/ 		pZ)$ on a two-dimensional compact metric 
space $X$         	with $n$-dimensional orbit space.}
 
This theorem was proved in [DW] with an error in one of the lemmas
(Lemma 15). This paper presents a corrected version of Lemma 15
and it is identical with [DW] in the rest.
	\endabstract \endtopmatter

\document

\head Introduction \endhead

	Interest in group actions with the property that the orbit space 
has greater dimension than the original space is motivated in part by P. 
A. Smith's generalization [S2] of Hilbert's Fifth Problem [H].  Hilbert 
asked whether a topological group $G$ the underlying space of which, 
$|G|$,  is a (finite- dimensional) manifold admits the structure of a 
Lie group, i.e., whether $|G|$ admits the structure of a real-analytic 
manifold on which group multiplication and inversion are analytic 
functions).  Hilbert's problem was solved affirmatively by von Neuman in 
1933 [vN] for compact groups and by Montgomery and Zippin in 1952 [MZ], 
cf. [G], for locally compact ones.  Smith's generalization, which asks 
whether a compact topological group acting effectively (i.e., so that 
each element that is not the identity moves some point) on a finite-
dimensional manifold must be a (subgroup of a) Lie group, is still 
unanswered.  It is now known as The Hilbert-Smith Conjecture. 	  
	The Hilbert-Smith Conjecture is known to be equivalent to the 
conjecture that a compact topological group $G$ acting effectively on a 
compact finite-dimensional manifold $M^n$ cannot have arbitrarily small 
subgroups; moreover, one may restrict attention to the case that $G$ is 
totally disconnected, so it is equivalent to the proposition that if $G$ 
is zero-dimensional, then it is finite, cf. [S2; Section 10].  Work of 
Newman [N] and Smith [S1] eliminates arbitrarily small torsion and 
reduces the problem to consideration of the $p$-adic integers 
$A_p=\varprojlim\{Z/p^nZ,\phi _n\}$, where $\phi _n(g)=g^p$ is the p-th 
power mapping, for $p$ prime [S2; Section 10].

	C.T. Yang [Y] (cf. [R]) showed that the orbit space 
$M^n/A_p$ of an effective $A_p$-action on an an $n$-manifold 
must have integral homological dimension dim$_Z M/A_p = n+2$, 
and that for any locally compact, Hausdorff, homologically n-dimensional 
space $X$ supporting an effective $A_p$ action, dim$_Z X/A_p \leq n+3$.
Actions of this sort are termed {\it dimension raising}.

  Since for locally compact metric spaces $Y$ of finite covering 
dimension dim $Y$, we have dim$_Z\, Y =\,$dim $Y$ [A] (cf. [Ku], [D]), 
Yang's Theorems imply that if $X$ is $n$-dimensional metrizable, then dim $X/A_p \leq 
n+3$ or dim $X/A_p = \infty$, and that dim $M/A_p = n+2$ or dim $M/A_p = 
\infty$ if $M$ is an $n$-manifold.  Examples of interesting dimension-raising actions of $A_p$ on 
non-manifold $X$ with dim $X/G =$ dim$X+1$ or dim $X + 2$ were constructed in [K] 
and [RW] (cf. [W1], [W2]), and an example of a free $A_p$-action on a 
2-dimensional cell-like set (apparently unpublished) was constructed by 
Bestvina and Edwards (cf. [Ga]).

	In contrast, nothing has been known until now about the 
dimension-raising possibilities for actions of torsion groups on compact 
metric spaces.  We show in this paper that the restrictions on the 
dimension of orbit spaces of actions by $A_p$ do not apply: dimension 
can be raised by an arbitrary amount if the acting group has enough 
torsion. Specifically, we prove the following theorems.  (Henceforth, we 
use $\Z_p$ to denote the integers modulo $p$ and $G$ to denote $\prod 
_{i=1}^{\infty} (\Z_p)_i$.) 
 \proclaim{Theorem A} For each integer $n \ge 3$ and each prime $p$ 
there is an action of $G$ on a compact, two-dimensional metric space $X$ 
such that dim $X/G = n$; moreover, dim $X \times X = 3$. 
\endproclaim 

\proclaim{Theorem B}  For each prime $p$ there is an action of $G = 
\prod _{i=1}^{\infty}(\Z_p)_i$ on a two-dimensional metric compactum with 
strongly infinite-dimensional orbit space. 
\endproclaim

	Theorems A and B are proved as a consequence of our third 
theorem. 
 \proclaim{Theorem C} For each compact metric space $Y$, and each prime 
$p$ $Y$ is the orbit space of an action of $G$ on a metric compactum $X$ 
with dim$_{\Z_p} X = 1$.
 \endproclaim 

The paper [DW] with Theorems A,B,C was published in 1997. In 2001
Bob Edwards pointed us out that there is a gap in our Lemma 15.
Since the gap was fixable and the lemma was not essentual in our construction,
we did not worry about it until the recent
preprint of Z. Yang [Ya] which fixes our lemma. 
Since it is rather our responsibility, 
we decided to publish corrections ourselves.

The problem with our Lemma 15 is that it is overstated. 
In fact, it is not really needed for the construction.
What was used in the proofs of
Theorems A,B,C is a much weaker statement. 

In this version of [DW] we present a weak version of Lemma 15 
which works also for $Z_{p^k}$ for $k>1$ as well.  This allows to extend
automatically Theorems A, B and C to the groups $G=\prod_{i=11}^{\infty}Z_{p^{k_i}}$ as 
in [Ya].
We don't do it, since we want to emphasize on the torsion
groups. Besides, we believe that by the nature of our
construction, the induced $A_p$-action in Theorems A,B,C, 
$A_p\subset\prod_{i=11}^{\infty}Z_{p^{i}}$, does not raise the dimension.

In this paper we also do revision of the proof of Lemma 16 which contained 
some inessential errors.
The rest of the paper remains intact.

\head Preliminaries \endhead

	Here we synopsize the dimension theory needed for our argument.  
Our basic references are [HW], [Ku] and [D].  From this point on, all 
spaces are locally compact and metric. 
 \definition{Definition}  A space $X$ has {\it covering dimension} dim 
$X \leq n$ provided there are arbitrarily fine open coverings $\Cal U$ 
of $X$ of order $\le n+1$, i.e., with the property that at most $n+1$ 
elements of $\Cal U$ have a non-void mutual intersection; otherwise dim 
$X > n$. We use the term $n$-{\it dimensional} in this sense. 
\enddefinition 
 \definition{Definition}  The {\bf Urysohn n-diameter} of a space $X$ is 
the infimum $d_n(X)$ of all positive numbers $\epsilon$ such that there 
is an $\epsilon$-map $f:X \to Y$ for some $n$-dimensional space $Y$.  
(An $\epsilon$-{\it map} is one with all point inverses of diameter less 
than or equal to $\epsilon$.) 
\enddefinition 
\proclaim{Proposition 1 } A compact metric space $X$ has Urysohn 
$n$-diameter zero if and only if dim $X \le n$. \endproclaim 
%
%
%

%
\proclaim{Theorem 2 (Hurewicz [HW])}  Let $f:X\rightarrow Y$ be a closed 
surjective mapping.  Then dim $X \le {\text dim}\; Y + {\text sup \; } 
\{\text{dim}\,  f^{-1}(y) \mid y \in Y\}$. 
\endproclaim 
\proclaim{Corollary 3} dim $X \times Y \leq {\text dim}\, X + {\text 
dim}\, Y$.
\endproclaim 
\proclaim{Corollary 4}  If $f:X @>>> Y$ is a continuous surjection of 
compact metric spaces with zero-dimensional point inverses, then dim $Y 
\ge \text{dim }\thickspace X$.
\endproclaim
%
%

\definition{Definition}  If $X$ is a locally compact metric space, its 
{\it cohomological dimension} dim$_G\, X$ with respect to an Abelian 
group $G$ is the supremum of the set of integers $n$ for which there 
exists a locally compact set $A \subseteq X$ with $H^n_c(A;G) \neq 0$. 
\enddefinition
	Here, cohomology with compact supports $H^*_c(A;G)$ is \v Cech 
cohomology of the one-point compactification of $A$ if $A$ is 
non-compact (and \v Cech cohomology of $A$ if $A$ is compact). 
 \proclaim{Theorem 5 (Cohen [C], cf. [Ku])}  For a locally compact 
metric space $X$ and Abelian group $G$ the following are equivalent:
 \roster 	
\item dim$_G\,X \le n$. 	
\item $H^{n+1}_c(U;G) = 0$ for all precompact open sets $U$ of $X$. 
\item the restriction homomorphism $j^{\ast}:H^n_c(B;G) \rightarrow 
H^n_c(A;G)$ is surjective for each compact pair $(A,B)$ in $X$.
 \endroster
 \endproclaim 
%
%
\proclaim{Theorem 6 (The Fundamental Theorem of (co-) Homological 
Dimension \linebreak
(Alexandroff [A], cf. [Ku], [Wa], [D]))}  Let $X$ be a locally 
compact metric space.  Then dim $_{Z}\,X \le \text{ dim}\, X$.  If dim 
$X < \infty$ then dim $_{Z}\, X =\text{dim}\, X$. 
\endproclaim 
Let $\Bbb Q$ denote the rational numbers. For a prime $p$, let $Z_{(p)}$ 
be the integers localized at $p$.  We also need 
$$Z_{p^{\infty}} = \varinjlim {\Z_p \rightarrow Z_{p^2} \rightarrow \dots 
\rightarrow Z_{p^n} \rightarrow \dots}$$

	Let $\sigma=\{\Bbb Q\}\bigcup _{p \; prime} \{\Z_p, Z_{(p)}, 
Z_{p^{\infty}}\}.$  For a group $G$, let $\sigma (G) \subseteq \sigma$ 
be the subcollection containing $\Bbb Q$ if and only if $G$ contains an 
element of infinite order, $\Z_p$ if and only if $G$ contains an element 
of order $p^k$ not divisible by $p$, $Z_{p^{\infty}}$ if and only if $G$ 
contains an element of order $p$, and $Z_{(p)}$ if and only if $G$ 
contains an element $a$ such that for all $n$, $p^na$ is not divisible 
by $p^{n+1}$ (writing group operations additively).  Then we have the 
following. 
\proclaim{Proposition 7 (Bockstein's Inequalities [Ku])}  For any 
locally compact space $X$ the following inequalities hold: \roster 
	\item dim$_{Z_{p^{\infty}}}\, X \le \text {dim}_{\Z_p}\, X$, 
	\item dim $_{\Z_p}\, X \le \text {dim}_{Z_{p^{\infty}}}\, X + 1$, 
	\item dim$_{\Bbb Q}\, X \le \text {dim}_{Z_{(p)}}\, X$, 
	\item dim$_{\Z_p}\, X \le \text {dim}_{Z_{(p)}}\, X$, 	\item 
dim$_{Z_{p^{\infty}}} \le \text {max}\, \{\text {dim}_{\Bbb Q}\, 
X,\thickspace \text {dim}_{Z_{(p)}}\, X-1\}$, 	\item dim$_{Z_{(p)}}\, X 
\le \text {max}\, \{\text {dim\, }_{\Bbb Q}\, X,\thickspace \text 
{dim}_{Z_{p^{\infty}}}\, X+1\}$. \endroster 
\endproclaim 
\proclaim{Theorem 8 (Realization Theorem [D,D1])} For every set of numbers
$\{n_G:\,G\in\sigma\}$ satisfying Bockstein's Inequalities there exists
a compactum $X$ with $dim_GX=n_G$ for all $G\in\sigma$.
\endproclaim
\proclaim{Theorem 9 (Bockstein's First Theorem [Ku])}  For any locally 
compact metric space $X$ and Abelian group $G$, dim$_G \, X = {\text 
max}_{H\in \sigma (G)}\,{\text dim}_H\, X$. \endproclaim %
\proclaim{Theorem 10 (Bockstein's Second Theorem [Ku])}  For $X$ and $Y$ 
locally compact metric spaces, the following hold: 	
\roster 		
\item dim$_{\Z_p}(X \times Y)=\text {dim}_{\Z_p}\,X+ \text {dim}_{\Z_p}\, 
Y$, 		
\item dim$_{\Bbb Q}\,(X\times Y) = \text {dim}_{\Bbb Q}\, X + \text 
{dim}_{\Bbb Q}\, Y$, 		
\item dim$_{Z_{p^{\infty}}}(X \times Y)= \text {max} \, \{\text 
{dim}_{Z_{p^{\infty}}}\, X + \text {dim}_{ Z_{p^{\infty}}}\, Y, 
\thickspace \text {dim}_{\Z_p}(X \times Y)-1\}$, 		
\item dim$_{Z_{(p)}}(X \times Y) = \text {dim}_{Z_{(p)}}\,X + \text 
{dim}_{Z_{(p)}}\,Y$, if dim$_{Z_{p^{\infty}}}\,X = \text 
{dim}_{Z_{(p)}}\,Y$ or if \linebreak dim$_{Z_{p^{\infty}}}\,Y = \text 
{dim}_{Z_{(p)}}\,X$, 		
\item dim$_{Z_{(p)}}(X \times Y) = \text {max }\{\text 
{dim}_{Z_{p^{\infty}}}(X \times Y)+1,\thickspace \text {dim}_{\Bbb Q} (X 
\times Y)\}$, if dim$_{Z_{p^{\infty}}}\,X < \text {dim}_{Z_{(p)}}\,X$ 
and dim$_{Z_{p^{\infty}}}\,Y < \text {dim}_{Z_{(p)}}\,Y$. 
	\endroster 
\endproclaim 

\proclaim{Theorem 11 (Countable Union Theorem [Ku]}
 If a compact space $X = \bigcup _{i=1}^{\infty}X_i$, where 
each $X_i$ is closed in $X$, then $dim_GX = \sup\{ dim_GX_i\}$.
\endproclaim

\head Lemmas \endhead

\proclaim{ Lemma 12}  If $f: X @>>> Y$ is a light continuous map of a compact 
metric space $X$ onto a finite dimensional $Y$, then for any Abelian 
group $G$, dim$_G\, X \le \text{dim}_G \,Y$. \endproclaim 
 \demo{Proof}  
Let dim $Y = n$.  For each Abelian group G there is a "Test space" 
$T_n(G)$ such that if dim $\,Z \le n$, then dim$_G \, Z = \text{dim} Z 
\times T_n(G) - n$.   By Corollary 2, dim $X \le n$ and dim $X \times 
T_n(G) \le \text{dim} Y \times T_n(G)$.  The result follows.
\enddemo 
\proclaim{Lemma 13}  If dim$_{\Z_p}\, X \le 1$ and dim$_{Z[\frac 
{1}{p}]}\, X \le 1$, then dim$_{Z_q}\, X \le 1$ for all primes $q$ and 
dim$_{\Bbb Q}\, X \le 1$. 
\endproclaim 
 \demo{Proof} First, $\sigma(Z[\frac {1}{p}])=\{\Bbb Q,\} \bigcup 
\{Z_{(q)}\mid q\neq p, \thickspace q\thickspace \text {prime}\}$.  By 
Theorem 5, dim$_{\Bbb Q}\, X \le 1$ and dim$_{Z_{(q)}}\, X \le 1$ for 
all primes $q \neq p$.  Proposition 2(4) now gives dim$_{Z_q}\, X \le 1$ 
for all primes $q$. 
\enddemo 
\proclaim{Lemma 14} The following are equivalent for a compact metric 
space $X$: 	
\roster 		
\item dim$_{F}\, X \le n$ for all fields $F$, 		
\item dim$_{F}\, X \le n$ for all fields $\Bbb \Z_p$ and for $\Bbb Q$, 
\item dim$_{Z}\, X\times X \le 2n+1$. 	
\endroster 
\endproclaim 
\demo{Proof}  That (1) implies (2) is trivial.  To see that (2) implies 
(1), note that if $F$ is a field of characteristic zero, $\sigma(F) = 
\{\Bbb Q\}$, while if $F$ has characteristic $p$, then $\sigma(F) = 
\{\Z_p, Z_{p^{\infty}}\}$ and by Proposition 2(1), 
dim$_{Z_{p^{\infty}}}\, X \le \text {dim}_{\Z_p}\, X.$

	To verify that (2) implies (3), note that $\sigma(Z) = \{\Bbb Q, 
\thickspace Z_{(p)} \text {for all primes}\thickspace p\}$, so applying 
Theorem 6 (1) and (2), Proposition 2 (6) and (5), and Theorem 6 (4) and 
(5) gives the result.

	Theorem 6 (1) and (2), together with Proposition 2(4) and the 
fact that $\sigma(Z) = \{\Bbb Q, \thickspace Z_{(p)} \text {for all 
primes}\thickspace p\}$ immediately show that (3) implies (2). 
\enddemo 

	We shall use $|L|$ to denote the underlying polyhedron of the 
simplicial complex $L$ when it is necessary to avoid confusion; 
otherwise, we shall use $L$ for both the complex and the polyhedron. 

\proclaim {Lemma 15}
For every finite connected simplicial complex $L$ with a
discrete group
$G$ simplicially acting on it there exists a
regular $Z_p^l$-cover
$p:M\to L$ such that
\roster
\item{} $p^*:H^1(L;Z_p)\to H^1(M;Z_p)$
is zero homomorphism;
\item{} If the fixed point set $L^G$
is nonempty and the action of $G$ on
$H_1(L;Z_p)$ is trivial, then the
action of $G$ on $M$ comutes with the deck
transformations.
\endroster

\endproclaim

\demo{Proof} For every complex $L$ the composition of the 

Hurewicz homomorphism in dimension one
and the tensor product with
$Z_p$ defines an epimorphism

$\xi_L:\pi_1(L)\to H_1(L;Z_p)=\oplus_{i=1}^lZ_p=Z_p^l=A$.

Its kernel $K=ker\xi_L$ defines
a regular $Z_p^l$-covering $p:M\to L$ of $L$. 

(1) Since $p_*(\pi_1(M))=K$
and $\xi_M$ is an epimorphism, it follows that

$p_*:H_1(M;Z_p)\to H_1(L;Z_p)$ is a zero homomorphism. Hence the dual
homomorphism $p^*$ is zero.

(2) Let $x_0\in L^G$ be a fixed point.
The space $M$ can be described as the space of classes $\{\gamma\}$ of paths
$\gamma:I\to L$ with $\gamma(0)=x_0$ where two such paths 
$\gamma$ and $\gamma'$ are equivalent
if the homotopy class $[\gamma\ast\bar\gamma']$ lies in $K$. 
For every $g\in G$ we define an action of $g$ on $M$ by the formula
$g\{\gamma\}=\{g\circ\gamma\}$. The action of $A$ on $M$ can be described 
as following. For every $a\in A$ take $[\alpha]\in\pi_1(L)$ such
that
$\xi_L([\alpha])=a$ and define $a\{\gamma\}=\{\alpha\ast\gamma\}$.

It is easy to check that this does not depend on choice of $\alpha$.

Since $\xi_L$ is a $G$-equivariant homomorphism and the action of $G$
on
$A$ is trivial by the hypotheses, we can replace $\alpha$ with $g\alpha$,

$g\in G$ in the above
description of the $A$-action:
$a\{\gamma\}=\{
\alpha\ast\gamma\}=\{g\alpha\ast\gamma\}$.
To complete the proof we need to
check
that the actions of $G$ and $A$ on $M$ commute.
Indeed, $(a\circ g)\{\gamma\}=
a\{g\gamma\}=\{g\alpha\ast g\gamma\}=
\{g(\alpha\ast\gamma)\}=(g\circ a)\{\gamma\}$ for every $g\in G$, $a\in A$,
and
$\{\gamma\}\in M$.
\qed
\enddemo

\proclaim {Lemma 16} Given a finite,
connected simplicial complex 
$L$ and a finite abelian group $A=Z_p$,
there are an integer $m$, a complex 
$\hat L$ and a simplicial action

of $A^m$ on $\hat L$ such that the orbit space $\hat L/A^m$ is
a 
subdivision of $L$ and the orbit map $f:{\hat L}\to L$
has the
following property for each simplex $\sigma \in L$:

$$  \text{The inclusion }\  f^{-1}(\partial \sigma)\to f^ {-1}(\sigma )\ \
\text{induces an isomorphism on } H_1(-;Z_p).  
\tag*$$

\endproclaim

\demo{Proof}  We make an induction on dimension $n$ of $L$.

If $n=0$ or $1$, we use $m=0$ and $f=$ id.
We will prove by induction
the following:  
The
restriction of $f^{-1}$ to the 1-skeleton $L^{(1)}$ defines an imbedding

$j:L^{(1)}\to \hat L$ which induces an isomorphism
$$ j^*:H^1(\hat L;Z_p)\to
H^1(L^{(1)};Z_p).
\tag**$$

Now assume that the Lemma holds true together
with (**) for complexes of dimension less than $n$.  
Let $L$ be an $n$-dimensional complex. We assume that an $A^k$-resolution
over the $n-1$-skeleton is already constructed: 
$h:\widehat{L^{(n-1)}}\to L^{(n-1)}$.
Let $\Delta$
be an $n$-simplex in $L$.   
We shall construct $f_{\Delta}:  \hat\Delta\to
\Delta$ in such a way that 
$f|_{\widehat{\partial\Delta}}=h$ where 
$\widehat{\partial\Delta}=h^{-1}(\partial\Delta)$.
Then the resolution $\hat L$ will be obtained by gluing $\hat\Delta$
to $\widehat{L^{(n-1)}}$ for all $n$-simplices $\Delta$.
Since
$\Delta^{(1)}$ is fixed under $A^k$-action, in view of (**) 
the action of
$A^k$ on $H_1(\widehat{\partial\Delta};Z_p)$ is trivial. 
By Lemma 15, there
exists a $A^k$-equivariant simplicial map 
$p :M\to\widehat{\partial\Delta}$ that induces zero homomorphism for 
1-dimensional mod $p$ cohomology.
Moreover, $A^l$ acts on $M$ with orbit space 
$\widehat{\partial\Delta}$ and the orbit map $p$, and $A^k\oplus A^l=A^{k+l}$ acts
on $M$ with orbit space $\partial
\Delta$ .  
Let $C_p$ denote the mapping cone of $p$.
Then we
get a natural  projection
$f_{\Delta}:C_p\to C_{id_{\partial \Delta^n}}=\Delta ^n$.  
By the definition this is the resolution 
$\widehat{\Delta}$.

Now to verify (*), we must show 
that the inclusion homomorphism
$$i^*:H^1(f^{-1}(\sigma );Z_p)\to 
H^1(f^{-1}(\partial \sigma );Z_p)$$ is an
isomorphism for 
each $\sigma \in L$.   The condition (*) holds
for every $k$-simplex, $k<n$, by induction assumption.
For an $n$-simplex $\Delta$  we show that
$i^*$ is an isomorphism.
This follows immediately from the cohomology exact
sequence 
of the pair $(f^{-1}(\Delta ^n),f^{-1}(\partial \Delta
^n))=
(C_p,\widehat{\partial \Delta ^n})$.  The quotient space

$C_p/\widehat{\partial \Delta ^n}=\Sigma M$, being a 
suspension of a
connected space, is simply connected.  This implies that $i^*$ is
a
monomorphism.  The boundary homomorphism 
$$\delta :H^1(\widehat{\partial
\Delta ^n};Z_p)\to 
H^2(C_p/\widehat{\partial \Delta ^n};Z_p)$$ coincides
with the 
suspension of the homomorphism $p^*$ and is hence zero.  
Thus,
$i^*$ is an epimorphism.

We note that $f^{-1}$ restricted
to $L^{(1)}$ coincides with the restriction of
$h^{-1}$. Hence 
$f^{-1}$ defines an imbedding $j$ of $L^{(1)}$.
In view of
the formula $(f^{-1}\mid)^*=i^*\circ(h^{-1}\mid)^*$ this 
imbedding induces
an isomorphism of 1-dimensional mod $p$
cohomology by the induction assumption
and the fact that $i^*$ is an isomorphism.
Thus, the property (**) holds for
every $n$-simplex $\Delta$ in $L$.
For a finite $n$-dimensional complex $L$ we apply
induction on
the number of $n$-simplices in $L$ and the Mayer-Vietoris exact
sequence
to a decomposition $L=N\cup\Delta$
to obtain the commutative
diagramm for mod $p$ cohomology:
$$
\CD
0 @>>> H^1(\hat L) @>>> H^1(\hat
N)\oplus H^1(\hat\Delta) @>>> 
H^1(\widehat{N\cap\Delta})\\
@VVV @Vj^*VV
@Vj_1^*VV @Vj_2^*VV\\
0 @>>> H^1(L^{(1)}) @>>> H^1(N^{(1)})\oplus
H^1(\Delta^{(1)}) @>>> 
H^1((N\cap\Delta)^{(1)})\\
\endCD
$$
in which $j_1^*$
and $j_2^*$ are isomorphisms. Then $j^*$ is an isomorphism
by the Five Lemma.

\qed
\enddemo

\head Proofs of the Theorems \endhead

\proclaim{Theorem C} For each compact metric space $Y$, and each prime 
$p$ $Y$ is the orbit space of a $G$-action on a metric compactum $X$ 
with dim$_{\Z_p} X = 1$. \endproclaim 
\demo{Proof}  First we address 
the case $Y = \Delta ^n$.  We apply Lemma 16, obtaining a resolution 
$q_1:(L_1,\sigma _1) \rightarrow (Y,\tau _1)$, where $L_1$ is a 
polyhedron with triangulation $\sigma _1$, $\tau _1$ is a triangulation 
of $Y$, and $q_1$ is a simplicial map that is the orbit mapping of a 
simplicial action of $G_{m_1}= \prod _{i=1}^{m_1}(\Z_p)_i$ on $L_1$ 
satisfying $\ast$.  Now take the barycentric subdivision $\tau _1'$ of 
$\tau _1$ and   repeat the process over $(Y,\tau _1')$, 
obtaining $q_2:(L_2,\sigma _2) \rightarrow (Y,\tau _2)$, a simplicial 
orbit map of an action of $G_{m_2}=\prod _{i=1}^{m_2}\Z_p$ on $L_2$ to a 
refinement $\tau _2$ of $\tau _1'$ satisfying $\ast$ with respect to 
$\tau _1'$.  Repeat this process infinitely often, obtaining for each 
$i=1,2,\dots$, $q_i:(L_i,\sigma _i)\rightarrow (Y,\tau _i)$, where $q_i$ 
is a simplicial orbit map of a simplicial action of $G_{m_i}$, $q_i$ 
satisfies $\ast$ with respect to $\tau _{i-1}'$.

	Now construct a sequence of simplicial pull-backs as follows:  
Let $\nu _1$ denote the pull-back to $|L_1|$ of the triangulation $\tau 
_2$, and let $p_2:(P_2,\upsilon _2) \rightarrow (L_1,\nu _1)$ be the 
pull-back of $q_2:(L_2,\sigma _2) \rightarrow (Y,\tau _2)$ over 
$q_1:(L_1,\nu _1) \rightarrow (Y,\tau _2)$, with $h_1:P_2 \rightarrow 
L_2$.  Let $k_2 = q_2 \circ h_2:P_2 \rightarrow Y$.  Then $p_2, h_2$
 and $k_2$ are simplicial and $G_{m_1 + m_2}=G_{m_1} \oplus 
G_{m_2}$ acts simplicially on $(P_2,\upsilon _2)$.

	Let $\nu _2$ denote the pull-back of $\tau _3$ by $k_2$ and form 
the simplicial pull-back $p_3:(P_3,\upsilon _3) \rightarrow (P_2,\nu 
_2)$ over $k_2$, with $h_3:P_3 \rightarrow L_3$ and $k_3 = q_3 \circ 
h_3$.  All is simplicial and there is a simplicial action of 
$G_{m_1+m_2+m_3}$ on $(P_3,\nu _3)$.  Repeat the process infinitely 
often, obtaining, for $i=3,\dots p_i:(P_i,\upsilon _i) \rightarrow 
(P_{i-1},\nu _{i-1})$ and $k_i:(P_i,\upsilon _i)\rightarrow (Y,\tau 
_{i+1})$ such that $k_{i-1} \circ p_i = k_i:|P_i| \rightarrow Y$.  Form 
the inverse limit $X = \underset\leftarrow \to{\text{lim}} \{|P_{i- 
1}|\leftarrow |P_i|: p_i\}$ together with the map $k = \underset 
\leftarrow \to{\text{lim}} k_i: X \rightarrow Y$.  Then $G=\prod 
_{i=1}^{\infty}$ acts on $X$ and $k:X \rightarrow Y$ is the orbit 
mapping.

	To verify that dim$_{\Z_p}\,X \le 1$, it suffices by Theorem 2(2) 
to show that for each compact set $A \subset X$ the restriction map 
$j^{\ast}:H^1(X; \Z_p) \rightarrow H^1(A;\Z_p)$ is surjective, which is 
equivalent to showing that each map $\phi :A \rightarrow K(\Z_p,1)$ 
extends to $X$.  We want to use the fact that $K(\Z_p,1)$ is an Absolute 
Neighborhood Retract to relate $\phi$ to maps of $|P_i|$ into 
$K(\Z_p,1)$. Consider the infinite mapping cylinder $M_{\infty}$ of the 
inverse system $\{|P_{i-1}| @<p_i<< |P_i|\}_{i=3}^{\infty}$, that is, 
let $M_i$ be the mapping cylinder of $p_i$ and, for each i, identify the 
copy of $|P_{i- 1}|$ in $M_i$ with the domain end of $M_{i-1}$.  Let 
$C_j = \bigcup _ {i=1}^jM_i$.  Then the collapses generate an inverse 
system $C_{j-1}\leftarrow C_j$ containing $|P_{j-1}|\leftarrow |P_j|$ 
with inverse limit $Z = M_{\infty} \bigcup X$.  Now $A \subset Z$ and 
$\phi :A \rightarrow K(\Z_p,1)$ extends to a map $\Phi:U \rightarrow 
K(\Z_p,1)$ where $U$ is an open neighborhood of $A$ in $Z$.  

	The collapses of the mapping cylinders $M_i$ are connected to 
the identities of the $M_i$ by one-parameter families of retractions 
which extend to give a one-parameter family of retractions $r_t:Z 
\rightarrow Z$ of $Z$, $0 \le t \le 1$ with $r_1$ the identity and $r_0$ 
the retraction onto $|L_1| = |P_1|$ and such that for $t_n = (n-1)/n$, 
$r_{t_n}$ is the projection onto $C_n$.  By compactness, there is an $n$ 
such that $r_t(A) \subset U$ for all $t \ge t_n$, so Borsuk's Homotopy 
Extension Theorem implies that $\phi$ extends to $X$ provided that $\Phi 
\circ r_n$ extends to $|P_n|$.  For $n$ sufficiently large, the union 
$T$ of all subcomplexes of $(P_n,\nu _n)$ of the form $p_n^{- 
1}(\sigma)$ for simplices $\sigma$ of $|P_{n-1}|$ that intersect 
$r_{t_n}(A)$ lies in $U$.  Let $\psi = \Phi |_T :T \rightarrow K(\Z_p,1)$.  
We wish to extend $\psi$ to $|P_n|$.  We do it by induction over the 
sets $k_n^{-1}(\tau _n(Y)^i)$, the preimages of the $i$-skeleta of $\tau 
_n$.

	For a simplex $\eta$ of $\nu _{n-1}'$, we have $k_{n-1}$ is 
one-to-one on $\eta$ and carries it to a simplex $\zeta$ of $\tau _{n-
1}'$, which is the domain triangulation of $f_n$.  Therefore, $\ast$ 
yields that the restriction homomorphism $j^{\ast}:H^1(q_n^{-
1}(\zeta);\Z_p) \rightarrow H^1(q_n^{- 1}(\partial \zeta);\Z_p)$ is an 
epimorphism.  So $j^{\ast}:H^1(p_n^{- 1}(\eta);\Z_p) \rightarrow 
H^1(p_n^{-1}(\partial \eta);\Z_p)$ is an epimorphism. Thus, $\psi$ extends to 
$|P_n|$, and dim$_{\Z_p}X \le 1.$

	Now we construct an example over the Hilbert cube as follows. 
For each $n$ consider the projection $I^n \rightarrow I^{n+1}$ and form 
the diagram 
$$ 
\CD 
I^n         @<<<        I^{n+1}      @<q_{n+1}''<< X_{n+1}\\ 
@Aq_nAA                 @Ap_n''AA                  @AAq_{n+1}''A\\ 
\tilde X_n @<p_{n+1}'<< \hat X_{n+1} @<p_{n+1}''<< \tilde X_{n+1} 
\endCD 
$$
 where $q_n$ is assumed by induction, $q_1$ being obtained from the 
finite- dimensional version of Theorem C just proved, $\hat X_{n+1}$ 
is pull-back, $X_{n+1}$ is obtained by applying the finite-dimensional 
version of Theorem C, and $\tilde X_{n+1}$ is obtained by pullback.  Let 
$q_{n+1} = q_{n+1}' \circ q_{n+1}'': \tilde X_{n+1} @>>> I^{n+1}$ and 
let $p_{n+1} = p_{n+1}' \circ p_{n+1}'':\tilde X_{n+1} @>>> \tilde 
X_n$.  By induction, we have an action of $G$ on $\tilde X_n$, with orbit 
mapping 
with orbit map $q_{n+1}'$, so we get an induced action of $G \times G 
\cong G$ on $\tilde X_{n+1}$ with orbit map $q_{n+1}$. By Lemma 1 
dim$_{\Z_p} \tilde X_{n+1} \le \text{dim}_{\Z_p} X_{n+1} \le 1$. 	Now form 
the inverse limits 
$$ 
\CD I^1    @<<<    I^2        @<<<    I^3        @<<< \dots\\
 @AAq_1A           @AAq_2A            @AAq_3A              @.\\                                                                                                              
\tilde X_1 @<<p_2< \tilde X_2 @<<p_3< \tilde X_3 @<<< \dots  
\endCD 
$$
	 We now have an action of the infinite product $G^{\infty}\cong 
G$ on $X$ with orbit map $q$.  Again, Lemma 1 ensures that dim$_{\Z_p} X 
\le 1$.  To complete the proof, let $Y$ be any compact metric space, 
embed $Y$ in $I^{\infty}$, and let $X = q^{-1}(Y)$.

\enddemo                                 

\proclaim{Theorem A}  For each integer $n \ge 2$ and each prime $p$ 
there is an action of $G=\prod _{i=1}^{\infty}\Bbb \Z_p$ on a compact, 
two-dimensional, metric space $X$ such that dim $X/G = n$; moreover, dim 
$X \times X\leq 3$.

\endproclaim

\demo{Proof}  Fix $n$ and $p$.  If $n = 2$, let $X_2$ be a compact 
metric  space with dim $X_2$ $=2$ such that dim $X_2 \times X_2 = 3$ [Ku] 
and use the trivial $G$-action.  Otherwise, let $Y_n$ be an $n$-
dimensional compact metric space such that dim$_{\Bbb Z[\frac{1}{p}]}Y_n 
\le 1$ (see Theorems 4 and 5). 
  By Theorem C there is a compact, metric $X_n$ with 
dim$_{\Bbb Z_p} X_n \le 1$ equipped with a $G-$action   and orbit map 
$q_n:X_n @>>> Y_n$.  Then $X_n$ is the desired space, as follows. 
	Firstly, dim $X_n \le n$ by Corollary 2.  Next, dim$_{\Bbb 
Z[\frac{1}{p}]} X_n \le 1$ by Lemma 1.  Now, dim$_{\Bbb Z_q} X_n \le 1$ 
for all primes $q$ and dim$_{\Bbb Q} X_n \le 1$ by Lemma 2.  Now 
 Lemma 3 shows that 
dim$_{\Bbb Z} X_n \times X_n \le 3$.  Thus, by Theorem 3,  
dim $X_n \times X_n \le 3$. 

\enddemo

\proclaim{Corollary A}	For each prime $p$ there is an action of 
$G=\prod _{i=1}^{\infty}\Bbb Z_p$ on a compact, two-dimensional, metric 
space $X_{\infty}$ such that dim $X_{\infty} = 2$ and dim 
$X_{\infty}\times X_{\infty}=3$. \endproclaim

\demo{Proof}  Fix $p$.  For each integer $n \ge 3$, let $Y_n$ and $X_n$ be as 
in Theorem A. Consider the one-point compactifications $X_{\infty}$ and
$Y_{\infty}$ of the disjoint unions $\cup X_n$ and $\cup Y_n$. Now
dim$Y_{\infty}=\infty$. Theorem 7 implies that dim$X_{\infty}\times X_{\infty}$.
Clearly $G$ acts on $X_{\infty}$ with $Y_{\infty}$ as the orbit space.

\enddemo
\head Strongly Infinite-Dimensional Orbit Spaces \endhead

\proclaim{Definitions}
\roster
 

\item A map $f:X @>>> I^n$ is {\it essential} provided there is no homotopy 
of pairs $F:(X,f^{-1}(\partial I^n))\times I @>>> (I^n,\partial I^n)$ from $f$
to a map $g: X @>>> \partial I^n)$; otherwise, $f$ is {\it inessential}.

\item A map $f:X @>>> Q = \Pi _{i=1}^{\infty}$ into the Hilbert cube is {\it 
essential} provided that $p \circ f :X @>>> I^n$ is essential for each 
coordinate projection $p:Q @>>> I^n$; otherwise it is {\it inessential}.
 
 \item A space $X$ is {\it strongly infinite-dimensional} provided that there 
exists an essential map $f:X @>>> Q$ of $X$ onto the Hilbert cube.

 \item The symbol $X\tau Y$ means that $Y \in AE(\{X\})$, i.e., each map $f:A 
@>>>Y$ defined on a closed subspace of $X$ extends over $X$.

\endroster

\endproclaim
\proclaim{Lemma 6} For every prime $p$ there exists a strongly 
infinite-dimensional compactum $Y$ with the cohomological dimension
$dim_{\Bbb Z[\frac{1}{p}]}Y=1$.
\endproclaim
The proof of Lemma 6 is the same as the proof of Theorem 3 [D2] with the only
difference that in Lemma 6 one should use the cohomology $H^*(\ ;\Bbb Z_p)$
instead of K-theory $K^*_{\Bbb C}(\ ;\Bbb Z_p)$.
\proclaim{Theorem $B'$} For each prime $p$ there is an action of $G = \prod 
_{i=1}^{\infty}({\Bbb Z} _p)_i$ on a compact space $X$ with $dim_{\Bbb Z}X\times X\leq 3$
and strongly infinite dimensional orbit space.
\endproclaim
\demo{Proof}
Apply Theorem C to a compactum $Y$ of Lemma 6 to obtain a desired compactum $X$
together with $G$-action.
\enddemo
In the rest of this section we are improving Theorem $B'$ by constructing a finite
dimensional $X$. We need the following:
\proclaim{Theorem 8 (Generalized Eilenberg Theorem [D1]}  Let $K$ and $L$ 
be countable CW complexes and $X$ compact and metric.  Suppose that $X \tau K \ast 
L$.  Then for each map $g:A @>>> K $ of a closed subset of $X$ into $K$ there 
is a compact $Y \subset X - A$ such that $Y \tau L$ and  $g$ extends over $X 
- Y$.

\endproclaim
\proclaim{Theorem B}  For each prime $p$ there is an action of $G = \prod 
_{i=1}^{\infty}({\Bbb Z} _p)_i$ on a two-dimensional metric compactum with 
strongly infinite-dimensional orbit space.
\endproclaim
\demo{Proof}  Let $Q= I^{\infty}= \prod _{i=1}^{\infty} I_i$ and $Q_m = 
I_m^{\infty} = \prod _{i=n}^{\infty}I_i$.  Let $I_m^n = \prod _{i=m}^nI_i$.  With $I = 
[0,1]$ and base point $0$, we have natural inclusions.  Let $q_m^n$ denote the 
projection to $I_m^n$, for $m \le n \le \infty$ and $n < \infty$. 
	For each n, let $\eta _n: Z_n \to I^n$ be a resolution as in Theorem C of 
$I^n$ as orbit space of a $G$-action, where dim$_{{\Bbb Z}_p}Z_n \le 1$.  
Denote by $G_n$ the copy of $G$ that acts on $Z_n$.  Let $\zeta _n = q^{n-1} 
\circ \eta _n: Z_n \to I^{n-1}$, where $q^{n-1}=q^{n-1}_1$.  Set $\tilde{Z} _2 
= Z_2$, and let $k_2:\tilde{Z} _2 \to Z_2$ be the identity.  Let
$\tilde{\eta} _2 = \eta _2 \circ k_2$.
	Define $\tilde{Z} _n$ inductively as the pull-back of $\tilde{Z} _{n-1}$ 
and $Z_n$ over $I^{n-1}$, with induced maps $q^{n-1}:\tilde{Z}_n \to \tilde{Z} 
_{n-1}$ over $q^{n-1}$ and $k_n$ over $\tilde{\eta} _{n-1}$ as in $*$:

$$
\CD
I^{n-1}                                  @<\zeta _n<<                     Z_n \\
@A\tilde{\eta} _{n-1} AA                                   @AAk_n A  \\
\tilde{Z}_{n-1}                     @<<\tilde{q}^{n-1}<          \tilde{Z}_n 
\endCD \tag{$*$}
$$

	Setting $\tilde{\eta} _n = \eta _n \circ k_n$ sets up the next pull-back. 
Note that $\tilde{Z}_n$ inherits a $G_2 \times \dots \times G_n$-action 
from the pull-back construction, and $k_n:\tilde{Z}_n \to Z_n$ is the orbit 
map of the $G_2 \times \dots \times G_{n-1}$-action on $\tilde{Z}_n$.  
Thus, as $\eta _n: Z_n \to I^n$ is the orbit map of the $G_n$-action,  
$\tilde{\eta}_n : \tilde{Z}_n \to I^n$ is the orbit map of the entire $G_2 
\times \dots G_n$-action.  Moreover, $\tilde{q}^n$ is equivariant.
	We now have the morphism $\tilde{\eta}_n$ between 
the inverse systems $\{I^n, q^n \}$ and $\{\tilde{Z} _n, \tilde{q} 
^n \}$.	
	We shall construct inductively an inverse system $\{N_n, p_n \}$ with 
$N_n \subset I^n$ and $p_n = q^{n-1}|N_n$ so that the lift $\{\tilde{N} _n, 
\tilde{p} _n \}$
with $\tilde{N}_n  = \tilde{\eta} _n ^{-1} (N_n)$ and $\tilde{p}_n = \tilde{q} 
_n| \tilde{N} _n$ has the property that $\tilde{N}_n$ has Urysohn 2-diameter
$d_2(\tilde{N}_n) < 1/n$ (i.e., admits a map into a two-dimensional complex 
with point inverses of diameter less than $1/n$).  This ensures from the 
definition of covering dimension that $X = \underset\gets\to{\text{lim}} 
\tilde{N} _n$ is 2-dimensional.  The $N_n$'s will
be chosen carefully to ensure that $Y =  \underset\gets\to{\text{lim}} N_n$ is 
strongly infinite-dimensional.  To achieve this, we construct a (non-
inverse) sequence $Y_n \subset I^n$.
	Let $Y_2 \subset I^2$ be the set consisting of the center of $I^2$.  Now 
dim $\tilde{\eta} _2^{-1}(Y_2) =0$, so there is a map $\theta_2: \tilde{\eta} 
_2^{-1}(Y_2) \to K^2$ for some 2-complex $K^2$ with point inverses of 
diameter less than $1/2$ (e.g., the barycentric map to the nerve of an open 
cover of order at most $3$ with stars of diameter at most $1$ defined from 
a partition of unity subordinate to the cover).  Since $K^2$ is an Absolute 
Neighborhood Retract, $\theta _2$ extends over a neighborhood $U_2$ of 
$\tilde{\eta} _2^{-1}(Y_2)$ with point inverses of diameter less than $1/2$.  
Let $N_2 \subset \text{int}(I^2)$ be a closed disc containing $Y_2$ in its 
interior and sufficiently small that $\tilde{N}_2\tilde{\eta} _2^{-1}(N_2) 
\subset U_2$.
	Let $M = M(\Z_p,1)$ and $L=M(Z[\frac{1}{p}])$ be  CW Moore spaces (i.e., 
$H_i(M;Z)=\Z_p$ if $i=1$ and is otherwise zero, etc.).  Let $\phi _2: I^2-Y^2 
\to K$ be a map of degree 1.
	Assume for the induction that we have $Y_{n-1} \subset \text{int}(N_{n-
1}$, an essential map $\phi _{n-1}:N_{n-1}-Y_{n-1} \to K$, and 
$\tilde{N}_{n-1}= \tilde{\eta} _{n-1}^{-1}(N_{n-1})$.  Let $A_n = I^n- 
\text{int}(N_{n-1}\times I_n$.  Now define $\psi _n = i \circ \phi _{n-1} 
\circ q^{n-1} : A \to K \ast L$, where $i: K \to K \ast L$ is the natural 
inclusion.  From the Generalized Eilenberg Theorem (Theorem 8), we get a 
compact $Y_n \subset I^N - A = \text{int}(N_{n-1}) \times I_n$ such that 
$\psi _n$ extends to $I^n - Y_n$ and $Y_n \tau L$.  
	Since $L = M(Z[\frac{1}{p}]) \simeq K(Z[\frac{1}{p}])$, 
dim$_{Z[\frac{1}{p}]}Y_n \le 1$ by Theorem 4.1.  By Lemma 2, dim 
$_{Z[\frac{1}{p}]}\eta _n^{-1}(Y_n) \le 1$, and since dim$_{\Z_p}Z_n \le 1$, 
so is $ \eta _n^{-1}(Y_n)$.  It follows as in the proof of Theorem A that 
dim$\eta _n^{-1}Y_n \le 2$.  The maps $k_n$ all have zero-dimensional point 
inverses, so by Corollary 2, dim$\tilde{Y}_n = \tilde{\eta}_n^{-1}Y_n \le 2$, 
also.  Thus, as for $Y_2$, there is a two-complex $K^2_n$ and a map $\theta 
_n: \tilde{Y}_n \to K^2_n$ with point inverses of diameter less than $1/n$.  
With $U_n \subset \tilde{Y}_n$ a neighborhood of $\tilde{Y}_n$ such that 
$\theta _n$ extends to $U_n$ with point inverses of diameter less than 
$1/n$, choose $N_n \subset \text{int}N_{n-1} \times I_n $ a closed 
neighborhood of $Y_n$ such that $\tilde{\eta}_n^{-1}N_n = \tilde{N}_n 
\subset U_n$.  This completes the inductive construction.
	 Now we have $\tilde{\eta}: X=\underset\gets\to{\text{lim}}\tilde{N}_n \to 
\underset\gets\to{\text{lim}} N_n = Y \subset Q$.  Moreover, dim$X \le 2$ and
$X$  is equipped with a $\prod _{n=2}^{\infty}G_n = G$-action for which $\eta$
is  the orbit map.  It remains to verify the strong infinite dimensionality of 
$Y$.
	To demonstrate strong infinite dimensionality, it suffices to exhibit an 
essential mapping $f: Y \to Q$.  Let $f = q_2^{\infty}|_Y : Y \to Q_2$.  Now 
from the construction of $Y$, we may define
$$\phi : Q - Y = \bigcup _{n=2}^{\infty} (Q - N_n \times Q_{n+1}) \to K = 
M(\Z_p,1)$$
by $\phi | Q - N_n \times Q_{n+1} = \phi _n \circ q^n$.  Let $S^1 = \partial 
I^2$.  Then $\phi$ is essential, and the inclusion $S^1 \times Q_3 \to Q - Y$ 
induces an epimorphism $H^1(Q - Y;\Z_p) \twoheadrightarrow H^1(S^1 \times 
Q_3;\Z_p) \cong \Z_p$.  Were $f$ inessential, there would be $n \geq 2$ such 
that $q_3^n \circ f :Y \to I_3^n$ is inessential, and hence an $m \ge n$ such 
that $q_3^n|_{N_m} : N_m \to I_3^n$ is inessential (by the Homotopy Extension 
Theorem).  But then $f_2 = q_2|_{N_m} :N_m \to I_3^m$ would be also 
inessential and would deform to a map $f_2': N_m \to \partial I_3^m = 
S^{m-3}$, which by a standard application of the Homotopy Extension 
Theorem may be assumed to agree with $f_2$ on $B = N_m \cap I^2 \times 
S^{m-3}$.  From $f_2'$ we could obtain a map $r = id \times f_2' : N_m \to I^2 
\times S^{m-3}$ that is the identity on $B$.  Now consider the following 
diagram (with $\Z_p$ coefficients), where $D=S^1\times I_3^m, E=N_2\times 
I_3^m, \text{and}\  F=N_2\times S^{m-3}$:
$$
\CD
H^1(I^n-N_n)   @<AD<< H^n(N_m,B)   @<\delta<<   H^{n-3}(B)  @<i*<< H^{n-3}N_m)\\
@A\alpha 'AA            @A\alpha AA                      @AA\beta A          
@AAA\\ H_1(D)            @<AD<< H^{m-2}(E,F) @<\gamma<< H_{m-3}(F) @<<<   0
\endCD
$$
	
	Now, $\alpha '$ is dual to the homomorphism $H^1(I^m - N_m;\Z_p) \to 
H^1(D;\Z_p)$, which is an epimorphism, so $\alpha '$, hence $\alpha$, is an 
monomorphism.  Therefore, $\alpha \circ \gamma (1) \neq 0$.  However, 
$\alpha \circ \gamma (1) = \delta \circ i* \circ r* (1) = 0$, a contratiction 
establishing essentiality of $f$, hence the strong infinite dimensionality of 
$Y$, and the Theorem is proved.
\enddemo
\ \ \ \ \ \ \ \ \ \ \ \ \ \ \ \ \ \ \ \  \ \ \ \ \ \ \ \ \ \ \ \ \ \ \ REFERENCES	

[A] P.S. Alexandroff, {\it Dimensiontheorie, ein Betrag zur Geometire der 
abgeschlossen Mengen}, Math. Ann. 106 (1932), 161-238.

[Br] G.E. Bredon, {\it Introduction to compact transformation groups},
Academic Press, New York, 1972.

[BRW] G.E. Bredon, F. Raymond, and R.F. Williams, {\it p-Adic groups of
transformations}, Trans. Amer. Math. Soc. 99 (1961), 488-498.

[D] A.N. Dranishnikov, {\it Homological dimension theory}, Russian Mathematical Surveys
43 (1988), 11-63.

[D1] A.N. Dranishnikov, {\it On the mapping intersection problem}, Pacific
Journal of Mathematics,  Vol. 173 No 2, 1996, 403-412.

[D2] A.N. Dranishnikov, {\it Generalized cohomological dimension of compact metric spaces},
   Tsukuba J. Math. 14 (1990) no. 2, 247-262.

[DW] A.N. Dranishnikov, {\it Compact group action that raise dimension to
infinity}, Topology Appl. 80 (1997), 101-114.

[G] A. Gleason, {\it On the structure of locally compact groups}, Duke Math. J. 18 (1951),
65-104.

[Ga] D. Garity, {\it Some remarks on the Hilbert-Smith Conjecture}, Notes taken after
R.D. Edwards.

[H] D. Hilbert, {\it Mathematische Probleme}, Nachr. Akad. Wiss. Gottingen 1900, 253-297.

[HW] W. Hurewicz and H. Wallman, {\it Dimension Theory}, Princeton Mathematical Series,
Vol. 4, Princeton University Press, 1941.

[K] A.N. Kolmogorov, {\it Uber offene Abbildungen}, Ann. of Math. 38 (1937), 488-498.

[Ku] V.I. Kuz'minov, {\it Homological dimension theory}, Russian Mathematical Surveys
23 (1968), 1-45.

[MZ1] D. Montgomery and L. Zippin, {\it Small subgroups of locally compact groups},
Ann. of Math. 56 (1952), 213-241.

[MZ2] D. Montgomery and L. Zippin, {\it Topological Transformation Groups}, Interscience
Tracts in Pure and Applied Mathematics, Vol. 1, Interscience Publishers, Inc. New York,
1955.

[N] M.H.A. Newman, {\it A theorem on periodic transformations of spaces} Quart. J. Math. 2
(1931), 1-8.

[vN] J. von Neumann, Die Einfuhrung analytischer Parameter in topologischen Gruppen,
Ann. of Math. 34 (1933), 170-190.

[R] F. Raymond, {\it Groups of transformations on manifolds} Proc. Amer. Math. Soc. 12
(1961), 1-7. 

[RW] F. Raymond and R.F. Williams, {\it Examples of p-adic transformation groups}, Ann. of
math. 78 (1963), 92-106.

[S1] P.A. Smith, {\it Transformations of finite period, III, Newman's Theorem}, Ann. of
Math. 42 (1941), 446-458.

[S2] P.A. Smith, {\it Periodic and nearly periodic transformations} in Lectures in Topology
(R. Wilder and W. Ayres, eds.), University of Michigan Press, Ann Arbor, Mich., 1941,
pp. 159-190.

[W1] R.F. Williams, {\it A useful functor and three famous examples in topology}, Trans.
Amer. Math. Soc. 106 (1963), 319-329.

[W2] R.W. Williams, {\it The construction of certain 0-dimensional transformation groups},
Trans. Amer. Math. Soc. 129 (1967), 140-156.

[Wa] J.J. Walsh,{\it Dimension, cohomological dimension, and cell-like mappings}
Lecture Notes in Math., 870, Springer, 1981, pp. 105-118.

[Y] C.T. Yang, {\it p-Adic transformation groups}, Mich. Math. J. 7 (1960), 201-218.

[Ya] Z. Yang, {\it A compact group action which raises dimension to infinity},
 arXiv: math.GT/0212159v1, 11 Dec. 2002.
\enddocument